{\obeyspaces\global\let =\ }

\font\ttten=cmtt10

\newdimen\outputBaseLineSkip
\newskip\beginOutputSkip
\newskip\endOutputSkip
\def\looserOutput#1{%
  \advance\beginOutputSkip by #1
  \advance\endOutputSkip by #1
}
\def\tighterOutput#1{%
  \advance\beginOutputSkip by -#1
  \advance\endOutputSkip by -#1
}

\def\beginOutput{%
    \par
    \penalty -150
    \penalty -150
    \begingroup
      \def\\{%
          \leavevmode
          \hss
          \endgraf
          \penalty 150
          }
      \ttten
      \parindent = 24pt
      \def\${\char`\$}
      \def\{{\char`\{}
      \def\}{\char`\}}
      \catcode`\_=\the\catcode`a
      \catcode`\^=\the\catcode`a
      \catcode`\#=\the\catcode`a
      \catcode`\~=\the\catcode`a
      \catcode`\&=\the\catcode`a
      \parskip=0pt
      \lineskip=0pt
      \obeyspaces
      }

\def\endOutput{%
    \endgroup
    \par
    \penalty -150
    \penalty -150
    \noindent}

\documentclass[letter]{amsart}
\usepackage{amsmath, amsfonts, amssymb, amsthm, mathrsfs}
\usepackage{hyperref}
\newcommand{\Z}{\mathbb{Z}}
\newcommand{\Q}{\mathbb{Q}}
\newcommand{\B}{\mathcal{B}}

\newcommand{\I}{\mathcal{I}}
\newcommand{\rk}{\operatorname{rank}}
\newtheorem{prop}{Proposition}
\newtheorem{theorem}[prop]{Theorem}
\theoremstyle{definition}
\newtheorem*{mydef}{Definition}
\everymath{\displaystyle}

\begin{document}
\title{\bf{Matroids: a Macaulay2 package}}

\author{Justin Chen}
\address{Department of Mathematics, University of California, Berkeley,
California, 94720 U.S.A}
\email{jchen@math.berkeley.edu}

\begin{abstract}
We give an overview of the {\tt{Macaulay2}} package {\tt{Matroids}}, which introduces functionality to create and compute with matroids into {\tt{Macaulay2}}. Examples highlighting the use of many functions in the package are provided, including applications of matroids to other areas.
\end{abstract}


\vspace*{-0.3in}
\maketitle

\noindent
\textsc{Introduction.} A matroid is a combinatorial object which abstracts the notions of (linear algebraic, graph-theoretic) independence. Since their introduction by Whitney \cite{Whi}, matroids have found diverse applications in combinatorics, graph theory, optimization, and algebraic geometry, in addition to being studied as interesting objects in their own right. 

For the reader already familiar with matroids, this package\footnote{\url{https://github.com/jchen419/Matroids-M2}} provides capabilities to form matroids from a matrix, graph, or ideal; convert between various representations of matroids; create and detect existence of minors; compute Tutte polynomials and Chow rings; as well as applications of matroids to polyhedral and algebraic geometry, commutative algebra, optimization, and even group theory. Each of these will in turn be illustrated with examples in this article. Virtually all notation and results mentioned below can be found in Oxley \cite{Ox}.

One striking feature of matroids is the multitude of distinct ways to define them. This variety of equivalent  -- or \textit{cryptomorphic} -- ways to characterize matroids is a great strength of matroid theory, and one of the reasons for its ubiquity. From the perspective of this package, the key definition is via bases:

\begin{mydef}
Let $E$ be a finite set, and $\varnothing \ne \B \subseteq 2^E$ a set of subsets of $E$. The pair $(E, \B)$ is a \textit{matroid} if for any $B_1, B_2 \in \B$ and $b_1 \in B_1 \setminus B_2$, there exists $b_2 \in B_2 \setminus B_1$ with $B_1 \setminus \{b_1\} \cup \{b_2\} \in \B$.
\end{mydef}

The set $E$ is called the \textit{ground set} of the matroid $M = (E, \B)$, and $\B$ is the set of \textit{bases} of $M$. All bases have the same cardinality, called the \textit{rank} of $M$. Any subset of a basis is an \textit{independent} set. A subset of $E$ that is not independent is \textit{dependent}. The minimal (with respect to inclusion) dependent sets are \textit{circuits}. It is easy to see that any of bases, independent sets, dependent sets, and circuits determines the others.

As any subset of an independent set is independent, the set of independent sets of a matroid forms a simplicial complex on $E$, called the \textit{independence complex} of $M$, denoted by $\Delta_M$. Via Stanley-Reisner theory, $\Delta_M$ corresponds to a squarefree monomial ideal $I_{\Delta_M} := \langle \prod_{i \in C} x_i \mid C \text{ circuit} \rangle$, inside a polynomial ring $k[x_i \mid i \in E]$ (since faces of $\Delta_M$ are independent sets, the minimal nonfaces are precisely the minimal dependent sets, i.e. circuits). We call $I_{\Delta_M}$ the (circuit) ideal of $M$: internally, many algorithms in this package work directly with this ideal, to exploit {\tt{Macaulay2}}'s facility with monomial ideals. \\

\noindent
\textsc{A first example.} The most basic way to create a matroid is by specifying the ground set and list of bases:

\beginOutput
i1 : needsPackage "Matroids";\\
i2 : M = matroid(\{a,b,c,d\},\{\{a,b\},\{a,c\}\})\\
o2 = a matroid of rank 2 on 4 elements\\
o2 : Matroid\\
\endOutput

\noindent
This creates a matroid of rank 2 on the ground set $\{a, b, c, d\}$ with 2 bases. We can peek at the matroid to see more of its internal structure:

\beginOutput
i3 : peek M\\
o3 = Matroid\{bases => \{set \{0, 1\}, set \{0, 2\}\}\}\\
             cache => CacheTable\{...1...\}\\
             groundSet => set \{0, 1, 2, 3\}\\
             rank => 2\\
\endOutput

\noindent
Two things should be noticed: first, {\tt{groundSet}} is a set of integers $\{0, \ldots, 3\}$ (instead of the given list $\{a, b, c, d\}$). Second, the bases consist of a list of subsets of {\tt{groundSet}}. This convention is by design: internally, the ground set is always identified with the set $\{0, \ldots, |E|-1\}$, and all sets associated to the structure of the matroid are subsets of the ground set. One should think of the integers in {\tt{groundSet}} as \textit{indices} of the actual elements, so $0$ is the index of the first element (in this case $a$), $1$ is the index of the second element, etc. 

The actual elements of the user-inputted ground set are not lost though: they have been cached in the {\tt{CacheTable}}, and can be accessed by using indices as subscripts on $M$, or all at once with an asterisk:

\beginOutput
i4 : (M_3, M_\{0,1\}, M_(set\{1,2\}), M_*)\\
o4 = (d, \{a, b\}, \{b, c\}, \{a, b, c, d\})\\
\endOutput

So far, no attempt has been made to check that $M$ is actually a matroid. We verify this now using the method {\tt{isWellDefined}} (which internally checks the circuit elimination axiom), and also give a non-example.

\beginOutput
i5 : (isWellDefined M, isWellDefined matroid(\{a,b,c,d\},\{\{a,b\},\{c,d\}\}))\\
o5 = (true, false)\\
\endOutput

We can obtain plenty of matroid-theoretic information for this example. Recall the following definitions:
\begin{mydef}
A \textit{loop} in $M$ is a 1-element circuit, and a \textit{coloop} in $M$ is an element contained in every basis. For $A \subseteq E$, the \textit{rank} of $A$ is the size of the largest independent subset of $A$, and the \textit{closure} of $A$ is $\overline{A} := \{x \in E \mid \rk(A) = \rk(A \cup \{x\}) \}$. A \textit{flat} of $M$ is a closed subset, i.e. $A = \overline{A}$. A \textit{hyperplane} of $M$ is a flat of rank equal to $\rk M - 1$.
\end{mydef}

\beginOutput
i6 : (rank M, rank(M, set\{0,3\}))\\
o6 = (2, 1)\\
i7 : (circuits M, independentSets(M, 1))\\
o7 = (\{set \{1, 2\}, set \{3\}\}, \{set \{0\}, set \{1\}, set \{2\}\})\\
i8 : (loops M, coloops M, closure(M, set\{2,3\}), hyperplanes M)\\
o8 = (\{3\}, \{0\}, set \{1, 2, 3\}, \{set \{0, 3\}, set \{1, 2, 3\}\})\\
i9 : flats M -- sorted by increasing size\\
o9 = \{set \{3\}, set \{0, 3\}, set \{1, 2, 3\}, set \{0, 1, 2, 3\}\}\\
i10 : fVector M -- number of flats of rank i, for 0 <= i <= rank M\\
o10 = \{1, 2, 1\}\\
\endOutput
\vspace{0.05cm}

\noindent
\textsc{Constructing types of matroids.} The simplest family of matroids is the family of \textit{uniform} matroids, where the set of bases equals all subsets of a fixed size:

\beginOutput
i11 : U = uniformMatroid(2,4); bases U\\
o12 = \{set \{0, 1\}, set \{0, 2\}, set \{1, 2\}, set \{0, 3\}, set \{1, 3\}, set \{2, 3\}\}\\
o12 : List\\
\endOutput

Another family of fundamental importance is the class of \textit{linear} matroids, which arise naturally from a matrix. The columns of the matrix form the ground set, and a set of column vectors is declared independent if they are linearly independent in the vector space spanned by the columns.

\beginOutput
i13 : A = matrix\{\{0,4,-1,6\},\{0,2/3,7,1\}\},; MA = matroid A; representationOf MA\\
o15 = | 0 4   -1 6 |\\
\      | 0 2/3 7  1 |
\endOutput

An abstract matroid $M$ is called \textit{representable} or \textit{realizable} over a field $k$ if $M$ is \textit{isomorphic} to a linear matroid over $k$, where an \textit{isomorphism} of matroids is a bijection between ground sets that induces a bijection on bases. We verify that the matroid $M$ we started with is isomorphic to $MA$, hence is representable over $\Q$:

\beginOutput
i16 : areIsomorphic(M, MA)\\
o16 = true\\
\endOutput

A matroid can also be constructed by specifying its circuit ideal, which we do for the same $M$ above. Here two matroids are considered equal if they have the same set of bases and same size ground sets; or equivalently, the identity permutation is an isomorphism between them.

\beginOutput
i17 : R = QQ[x,y,z,w]; MI = matroid ideal(y*z, w)\\
o18 = a matroid of rank 2 on 4 elements\\
i19 : M == MI\\
o19 = true\\
\endOutput

An important class of representable matroids (over any field) is the class of \textit{graphic} matroids, derived from a graph. If $G$ is a(n undirected) graph, then the graphic matroid $M(G)$ has ground set equal to the edge set of $G$, and circuits given by cycles in $G$, including loops and parallel edges. 

\beginOutput
i20 : K5 = completeGraph 5; M5 = matroid K5\\
o21 = a matroid of rank 4 on 10 elements\\
i22 : #bases M5 -- == n^(n-2) for M(K_n), by Cayley's theorem\\
o22 = 125\\
\endOutput

In this package, the graphic matroid is created by specifying circuits. This can be done for an abstract matroid as well, using the optional argument {\tt{EntryMode => "circuits"}} in the constructor function. Regardless of the value of {\tt{EntryMode}}, the bases are automatically computed upon creation. We recreate the matroid $M$ from before, by specifying its circuits (note the similarity with specifying the circuit ideal): 

\beginOutput
i23 : M == matroid(\{a,b,c,d\},\{\{b,c\},\{d\}\}, EntryMode => "circuits")\\
o23 = true\\
\endOutput

Certain common matroids are close to uniform, in the sense that relatively few subsets of size $\rk M$ are dependent, so the set of \textit{nonbases} ($=$ dependent sets of size $\rk M$) can also be specified: 

\beginOutput
i24 : nb = \{\{0,1,2\},\{0,4,5\},\{0,3,6\},\{1,3,5\},\{1,4,6\},\{2,3,4\},\{2,5,6\}\}/set;\\
i25 : F7 = matroid(toList(0..6), nb, EntryMode => "nonbases")\\
o25 : a matroid of rank 3 on 7 elements\\
i26 : (#bases F7, #circuits F7)\\
o26 = (28, 14)\\
\endOutput

A few specific matroids of theoretical importance are also built-in: currently $F_7, F_7^{-}, V_8, V_8^{+}, AG(3, 2), R_{10}$, and the Pappus and non-Pappus matroids. A library of all matroids on up to 8 elements is included as well:

\beginOutput
i27 : F7 == specificMatroids "fano"\\
o27 = true\\
i28 : L5 = allMatroids 5 -- non-isomorphic matroids on 5 elements\\
o28 = \{a matroid of rank 0 on 5 elements, a matroid of rank 1 on 5 elements, ... \\
i29 : (#L5, #flatten apply(6, i -> allMatroids i))\\
o29 = (38, 70)
\endOutput

One can also construct a new matroid from smaller ones by taking \textit{direct sums}: if $M_1 = (E_1, \B_1), M_2 = (E_2, \B_2)$ are matroids, then their direct sum is $M_1 \oplus M_2 := (E_1 \sqcup E_2, \{B_1 \sqcup B_2 \mid B_1 \in \B_1, B_2 \in \B_2\})$. A matroid that cannot be written as a direct sum of nonempty matroids is called \textit{connected}. Every matroid is a direct sum of connected matroids, its \textit{connected components}, which are unique up to rearrangement:

\beginOutput
i30 : S = U ++ matroid completeGraph 3\\
o30 = a matroid of rank 4 on 7 elements\\
i31 : C = components S\\
o31 = \{a matroid of rank 2 on 4 elements, a matroid of rank 2 on 3 elements\}\\
i32 : S == C#0 ++ C#1 and C#0 == U and C#1 == matroid completeGraph 3\\
o32 = true\\
\endOutput
\vspace{0.05cm}

\noindent
\textsc{Duality and minors.} One of the most important features of matroid theory is the existence of a duality. It is straightforward to check that if $M = (E, \B)$ is a matroid, then $\{E \setminus B \mid B \in \B\}$ is the set of bases of a matroid on $E$, called the \textit{dual matroid} of $M$, denoted by $M^*$. 

\beginOutput
i33 : D = dual M; (bases M, bases D)\\
o34 = (\{set \{0, 1\}, set \{0, 2\}\}, \{set \{2, 3\}, set \{1, 3\}\})\\
i35 : M == dual D\\
o35 = true\\
\endOutput

Virtually any matroid-theoretic property or operation can be enriched by considering its dual version -- for instance, loops of $M^*$ are coloops of $M$, and circuits of $M^*$ are complements of hyperplanes of $M$ (this is in fact how the method {\tt{hyperplanes}} works). Another operation is deletion, which dualizes to contraction:

\begin{mydef} Let $M = (E, \B)$ be a matroid, and $S \subseteq E$. The \textit{restriction} of $M$ to $S$, denoted $M|_S$, is the matroid on $S$ with bases $\{B \cap S \mid B \in \B, |B \cap S| = \rk S\}$. The \textit{deletion} of $S$, denoted $M \setminus S$, is the restriction of $M$ to $E \setminus S$. The \textit{contraction} of $M$ by $S$, denoted $M/S$, is defined as $(M^* \setminus S)^*$.
\end{mydef}

\beginOutput
i36 : N1 = M \symbol{`\\} set\{3\}; (N1_*, bases N1)\\
o37 = (\{a, b, c\}, \{set \{0, 1\}, set \{0, 2\}\})\\
i38 : N2 = M / set\{1\}; (N2_*, bases N2)\\
o39 = (\{a, c, d\}, \{set \{0\}\})\\
\endOutput

A \textit{minor} of $M$ is any matroid which can be obtained from $M$ by a sequence of deletions and contractions. It is a fact that any minor of $M$ is of the form $(M / X) \setminus Y$ for disjoint subsets $X, Y \subseteq E$. 

\beginOutput
i40 : minorM5 = minor(M5, set\{9\}, set\{3,5,8\}) -- contracts \{9\}, then deletes \{3,5,8\}\\
o40 = a matroid of rank 3 on 6 elements\\
i41 : (minorM5_*, #bases minorM5)\\
o41 = (\{set \{0, 1\}, set \{0, 2\}, set \{0, 3\}, set \{1, 2\}, set \{1, 4\}, set \{2, 3\}\}, 16)\\
\endOutput

Minors can be used to describe many important classes of matroids. For example, a class $\mathcal{M}$ of matroids is said to be \textit{minor-closed} if every minor of a matroid in $\mathcal{M}$ is again in $\mathcal{M}$. The classes of uniform, $k$-representable (for any field $k$), and graphic matroids are all minor-closed. Various classes of matroids can be characterized by their \textit{forbidden} or \textit{excluded} minors: namely the matroids not in the class, but with every proper minor in the class. 

\begin{theorem}[Tutte~\cite{Tu1, Tu2}] Let $M$ be a matroid.

i) $M$ is \textit{binary} ($=$ representable over $\mathbb{F}_2$) iff $M$ has no $U_{2,4}$ minor (i.e. no minor of $M$ is isomorphic to $U_{2,4}$).

ii) $M$ is regular ($=$ representable over any field) iff $M$ has no $U_{2,4}$, $F_7$, or $F_7^*$ minor.

iii) $M$ is graphic iff $M$ has no $U_{2,4}$, $F_7$, $F_7^*$, $M(K_5)^*$, or $M(K_{3,3})^*$ minor.

\noindent
Here $U_{2,4}$ is the uniform matroid of rank 2 on 4 elements, and $F_7$ is the Fano matroid.
\end{theorem}

We illustrate this by verifying that $M(K_5)$ is regular (alternatively, note that for any graph $G$, the signed incidence matrix of any orientation of $G$ represents $M(G)$ over any field):

\beginOutput
i42 : any(\{U, F7, dual F7\}, forbidden -> hasMinor(M5, forbidden))\\
o42 = false\\
\endOutput

Every minor of $M$ is in fact of the form $(M / I) \setminus I^*$, where $I, I^*$ are disjoint, $I$ is independent, and $I^*$ is \textit{coindependent} ($=$ independent in $M^*$). Such a minor has rank equal to that of $M/I$, which is equal to $\rk M - |I|$. Thus checking existence of a minor $N$ in $M$ can be realized as a two-step process, where the first step contracts independent sets of $M$ of a fixed size down to the rank of $N$, and the second step deletes coindependent sets down to the size of $N$.


\beginOutput
i43 : M4 = matroid completeGraph 4; hasMinor(M5, M4)\\
Contract set \{9\}, delete set \{3, 5, 8\}\\
o44 = true\\
i45 : minorM5 == M4\\
o45 = true\\
\endOutput

Finally, the \textit{Tutte polynomial} $T_M(x,y)$ of a matroid is an invariant which is universal with respect to satisfying a \textit{deletion-contraction recurrence}. It is a bivariate polynomial over $\Z$ which can be defined by the relation
\[
\hspace{5cm} T_M(x,y) = T_{M \setminus e}(x,y) + T_{M/e}(x,y), \qquad e \in E \text{ not a loop or coloop }
\]
with the initial condition $T_M(x,y) = x^ay^b$ if $M$ consists of $a$ coloops and $b$ loops. Any numerical invariant of matroids which satisfies a (weighted) deletion-contraction recurrence is an evaluation of the Tutte polynomial, up to a scale factor. For instance, the number of bases is equal to $T_M(1,1)$:

\beginOutput
i46 : tuttePolynomial M5\\
\       6     5    4       3      4     3      2         2      3      2              2\\
o46 = y  + 4y  + x  + 5x*y  + 10y  + 6x  + 10x y + 15x*y  + 15y  + 11x  + 20x*y + 15y  $\cdot\cdot\cdot$\\
i47 : tutteEvaluate(M5, 1, 1)\\
o47 = 125\\
\endOutput

For graphic matroids, the Tutte polynomial contains a wealth of information about the graph; e.g. the Tutte polynomial specializes to the chromatic polynomial. Even evaluations at specific points contain nontrivial information: e.g. $T_{M(G)}(2,1)$ counts the number of spanning forests in $G$, and $T_{M(G)}(2,0)$ counts the number of acyclic orientations of $G$.

\beginOutput
i48 : (tutteEvaluate(M5, 2, 1), tutteEvaluate(M5, 2, 0), factor chromaticPolynomial K5)\\
o48 = (291, 120, (x)(x - 4)(x - 3)(x - 2)(x - 1))\\
\endOutput
\vspace{0.05cm}

\noindent
\textsc{Connections.} We now present some connections of matroids to other areas of mathematics. First, polyhedral geometry: 
let $M = ([n], \B )$ be a matroid on $\{1, \ldots, n\}$. In Euclidean space $\mathbb{R}^n$ with standard basis $\{e_1, \ldots, e_n\}$, define the matroid polytope $P_M$ by taking the convex hull of the indicator vectors of the bases of $M$: 
\[
P_M := \operatorname{conv} \left( \sum_{i \in B} e_i \mid B \in \B \right)
\]

\noindent
The matroid polytope can be created as follows:

\beginOutput
i49 : needsPackage "Polyhedra"; P = convexHull basisIndicatorMatrix M4\\
o50 = \{ambient dimension => 6           \}\\
\       dimension of lineality space => 0\\
\       dimension of polyhedron => 5\\
\       number of facets => 16\\
\       number of rays => 0\\
\       number of vertices => 16\\
o50 : Polyhedron\\
\endOutput

\noindent
A theorem of Gelfand, Goresky, MacPherson, and Serganova \cite{GGMS} classifies the subsets $\B \subseteq 2^{[n]}$ which are the bases of a matroid on $[n]$ in terms of the polytope $P_M$. 

Next is optimization: let $E$ be a finite set, and $\I \subseteq 2^E$ a set of subsets that is downward closed: if $X \in \I$ and $Y \subseteq X$, then $Y \in \I$. Let $w$ be a weight function on $E$, i.e. a function $w : E \to \mathbb{R}$, extended to $w : 2^E \to \mathbb{R}$ by setting $w(X) := \sum_{x \in X} w(x)$. Consider the optimization problem (*) of finding a maximal member of $\I$ of maximum weight, with respect to $w$. One attempt to solve (*) is to apply the greedy algorithm: namely, after having already selected elements $\{x_1, \ldots, x_i\}$, choose an element $x_{i+1} \in E$ of maximum weight such that $\{x_1, \ldots, x_i, x_{i+1}\} \in \I$, and repeat. It turns out that the greedy algorithm will always work iff $\I$ is the set of independent sets of a matroid:

\begin{theorem} \cite{Bor}
Let $E$ be a finite set, and $\I \subseteq 2^E$. Then $\I$ is the set of independent sets of a matroid on $E$ iff $\I$ is downward closed and for all weight functions $w : E \to \mathbb{R}$, the greedy algorithm successfully solves (*).
\end{theorem}

\noindent
A solution to (*) provided by the greedy algorithm can be obtained using the method {\tt{maxWeightBasis}} (the weight function $w$ is specified by its list of values on $E$):

\beginOutput
i51 : w = \{0, log(2), 4/3, 1, -4, 2, pi_RR\}; maxWeightBasis(F7, w)\\
o52 = set \{3, 5, 6\}\\
\endOutput

Another application to optimization comes from the operation of \textit{matroid union}: if $M_1, M_2$ are matroids with independent sets $\mathcal{I}_1, \mathcal{I}_2$, then the independent sets of the union are of the form $I_1 \cup I_2$, where $I_1 \in \mathcal{I}_1, I_2 \in \mathcal{I}_2$ (and thus coincides with the direct sum if the ground sets are disjoint). 

\beginOutput
i53 : matroid(\{a,b,c,d\}, \{\{a\},\{b\},\{c\}\}) + matroid(\{a,b,c,d\}, \{\{b\},\{c\},\{d\}\}) == U\\
o53 : true\\
i54 : F7 + F7 == uniformMatroid(6, 7)\\
o54 : true\\
\endOutput

\noindent
Matroid union is an important operation in combinatorial optimization, and is closely related to transversal and matching problems: a matroid is \textit{transversal} iff it is a union of rank 1 matroids, and \textit{gammoids} (a class of matroids defined from vertex paths in directed graphs) are the minor-closure of the transversal matroids.

One can also find connections to group theory via the method {\tt{getIsos}}, which computes all isomorphisms between two matroids. Many interesting groups can be realized as automorphism groups of small matroids:

\beginOutput
i55 : aut = getIsos(F7, F7)\\
o55 : \{\{0, 1, 2, 3, 4, 5, 6\}, \{1, 0, 2, 3, 4, 6, 5\}, \{0, 2, 1, 3, 5, 4, 6\}, \{2, 0, 1, ...\\
i56 : #aut\\
o56 : 168\\
\endOutput

\noindent
The above output is an explicit permutation representation of $\operatorname{Aut}(\mathbb{P}^2_{\mathbb{F}_2}) = \operatorname{PGL}(3, \mathbb{F}_2)$ as a subgroup of $S_7$. For a larger example, the automorphism group of the Steiner system $S(5, 6, 12)$ is the Mathieu group $M_{12}$, a sporadic simple group of order $95040 = 2^6 \cdot 3^3 \cdot 5 \cdot 11$. This in turn is also equal to the automorphism group of the realizable matroid associated to a particular $6 \times 12$ matrix over $\mathbb{F}_3$ (\cite{Ox}, p. 367), and a high-performance computing cluster took just under $2$ hours to compute the entire permutation representation of this group inside $S_{12}$. 

For an application to commutative algebra: matroids are closely related to the Cohen-Macaulay property, for symbolic powers of squarefree monomial ideals. Indeed, a theorem of Terai-Trung \cite{TT} states that if $I$ is a squarefree monomial ideal, then $I$ is the circuit ideal of a matroid iff every symbolic power $I^{(n)}$ is Cohen-Macaulay, for $n \ge 1$ (in fact, this is equivalent to requiring just $I^{(3)}$ to be Cohen-Macaulay). As one can quickly check whether an ideal is the ideal of a matroid, this can give a quick proof that a particular symbolic power is Cohen-Macaulay:

\beginOutput
i57 : M6 = matroid completeGraph 6; L = (irreducibleDecomposition ideal M6)/(P -> P^3); \\
i59 : try ( alarm 10; I3 = intersect L; ) -- doesn't finish in 10 seconds\\
i60 : time isWellDefined M6\\
\     -- used 0.359306 seconds\\
o60 : true
\endOutput

Last but not least is algebraic geometry; in particular the emerging field of combinatorial Hodge theory. For a matroid $M$ on ground set $E$ with no loops, one can define a Chow ring associated to $M$: for a field $k$, set
\begin{align*}
R := k[x_F \mid F \text{ proper, nonempty flat}&]/(I_1 + I_2), \\
I_1 := \left(\sum_{i_1 \in F} x_F - \sum_{i_2 \in F} x_F \; \Big| \; i_1, i_2 \in E \text{ distinct} \right)\!, \; I_2 := &\left(x_Fx_{F'} \mid F, F' \text{ incomparable} \right)
\end{align*}

\noindent
where $F, F'$ run over all nonempty proper flats of $M$. Then $R$ is a standard graded Artinian $k$-algebra of Castelnuovo-Mumford regularity $r := \rk M - 1$. A result of Adiprasito, Katz, and Huh \cite{AKH} states that $R$ is a Poincare duality algebra (in particular, is Gorenstein) and has the strong Lefschetz property: for general $l \in R_1$ and $j \le r/2$, multiplication by $l^{r - 2j}$ is an isomorphism $R_j \xrightarrow{\sim} R_{r - j}$. We illustrate the Gorenstein property for the Vamos matroid (which is a smallest matroid not realizable over any field), and conclude by computing the dual socle generator or \textit{volume polynomial} (which generates the Macaulay inverse system of $R$) for $M(K_4)$:

\beginOutput
i61 : V = specificMatroids("vamos"); (rank V, #V.groundSet, #bases V, #flats V)\\
o62 = (4, 8, 65, 79)\\
i63 : I = idealChowRing V; apply(0..<rank V, i -> hilbertFunction(i, I))\\
o63 : Ideal of QQ[x   , x   , x   , x   , x   , x   , x   , x   , x      , x      , ...      \\
\                   \{7\}   \{6\}   \{5\}   \{4\}   \{3\}   \{0\}   \{2\}   \{1\}   \{6, 7\}   \{5, 7\}\\
o64 = (1, 70, 70, 1)\\
i65 : cogeneratorChowRing M4\\
        2       2       2       2       2       2                                   \\
o65 = 2t    + 2t    + 2t    + 2t    + 2t    + 2t    - 2t   t       - 2t   t       + ...\\
        \{5\}     \{4\}     \{3\}     \{2\}     \{1\}     \{0\}     \{5\} \{0, 5\}     \{0\} \{0, 5\}\\
\endOutput

\vspace{0.5cm}
\noindent
\textsc{Acknowledgements.} This project was partially supported by NSF grant DMS-1001867. The author is grateful to June Huh for explaining the connection to the Chow ring of a matroid. The author would also like to thank David Eisenbud and Daniel Grayson for advice with Macaulay2, Joe Kileel for enlightening discussions, Aaron Dall for software testing, Chris Eur for suggesting many valuable improvements, and the referee for helpful comments.

\end{document}